\documentclass{tran-l}

 \newtheorem{thm}{Theorem}[section]
 \newtheorem{cor}[thm]{Corollary}
 \newtheorem{lem}[thm]{Lemma}
 
 \theoremstyle{definition}
 \newtheorem{defn}[thm]{Definition}
 \newtheorem{defns}[thm]{Definitions}
 \theoremstyle{remark}
 \newtheorem{rem}[thm]{Remark}
 \newtheorem{rems}[thm]{Remarks}
 \newtheorem{exam}[thm]{Example}
 \newtheorem{exams}[thm]{Examples}
 \newtheorem{conj}[thm]{Conjecture}
 \newtheorem{que}[thm]{Question}
 \numberwithin{equation}{section}
 \usepackage[all]{xy}

\begin{document}

\title{Exchange maps of cluster algebras}
\date{}

\author{Ibrahim Saleh}

\maketitle
\begin{abstract}  Every two labeled seeds  in a field of fractions $\mathcal{F}$ together with a permutation give rise to an automorphism of $\mathcal{F}$ called an exchange map. We provide equivalent conditions for exchange maps to be cluster isomorphisms of the corresponding cluster algebras. The  conditions are given in terms of an action of the quiver automorphisms on the set of seeds.
\end{abstract}
\maketitle

\section{Introduction}

Cluster algebras were introduced by S. Fomin and A. Zelevinsky in [5, 6, 3, 7]. The original motivation was to find an algebraic combinatorial framework to study canonical basis and total positivity. One of the unique characterizations of cluster algebras is the way their generators (the cluster variables) are related. The  cluster variables are grouped in overlapping sets called clusters and each cluster forms a (commutative) free generators set of an ambient field. Attached with each cluster is a valued quiver and the cluster together with the valued quiver form a pair called a seed. A new seed is produced from an existing one using an operation called mutation, which defines an equivalence relation on the set of all seeds. The equivalence classes are called  mutation classes and each mutation class characterizes a cluster algebra.

In an effort to explore the automorphisms that preserve the mutation classes of coefficient-free skew-symmetric cluster algebras, the authors in [2] introduced and studied \emph{cluster automorphisms}, which are $\mathbb{Z}$-automorphisms of cluster algebras that send a cluster to another  and commute with  mutation. It was proved in [2, Corollary 2.7] that the cluster automorphisms are exactly the $\mathbb{Z}$-algebra automorphisms that map each cluster to a cluster. Inspired by this result and by the strong isomorphisms introduced in [6],  we mean by a \emph{cluster isomorphism} a $\mathbb{Z}$-algebra isomorphism that induces a bijection between the two sets of clusters of the two cluster algebras.\\
Any two seeds $(X, Q)$ and $(Y, Q')$ of rank $n$ in a field of fractions $\mathcal{F}$ and  a permutation $\sigma\in \mathfrak{S}_{n}$, define an automorphism over $\mathcal{F}$, induced by  sending the cluster variable $x_{i}$ in $X$ to the cluster variable $y_{\sigma (i)}$ in $Y$ for every $i$.  Such filed automorphism is  called  \emph{exchange map}, (Definition 3.10). The aim of this work is to find equivalent conditions on exchange maps to be cluster isomorphisms.\\ Quiver automorphisms act on the set of all seeds which gives rise to a relation between seeds called  \emph{$\sigma$-similarity} (Definition 3.5).
 \begin{thm} Two  cluster algebras of the same rank $n$ are cluster isomorphic if and only if there exists a permutation $\sigma\in \mathfrak{S}_{n}$ such that the  cluster algebras  contain two $\sigma$-similar seeds.
\end{thm}

   \emph{Positive cluster algebras}  are the  cluster algebras that satisfy the positivity conjecture [5]. The class of positive cluster algebras includes  all skew-symmetric cluster algebras [11], acyclic cluster algebras [10] and cluster algebras arising from triangulations of surfaces [12].

\begin{thm} An exchange map from a positive cluster algebra $\mathcal{A}$ to a positive cluster algebra $\mathcal{B}$  is a cluster isomorphism if and only if it sends every cluster variable in $\mathcal{A}$ to a cluster variable in $\mathcal{B}$.
\end{thm}
The article is organized as follows. Section 2 is devoted to cluster algebras associated with valued quivers. In the first three subsections of section 3, we introduce an action of quiver automorphisms on the set of all seeds, define the $\sigma$-similarity relation and provide some equivalent conditions for two seeds to be  $\sigma$-similar. In  subsection 3.4 we prove the main results (Theorems 3.14-3.15).
We finish section 3, by providing presentations of  groups of cluster automorphisms of some cluster algebras of types $\mathbb{B}_{2}$ and $\mathbb{G}_{2}$.

 Throughout the paper,  $K$ is a field with zero characteristic  and $\mathcal{F}=K(t_{1}, t_{2}\ldots t_{n})$ is the field of rational functions in $n$ independent (commutative) variables over $K$. Let $Aut_{K}(\mathcal{F})$ denote the automorphism group of $\mathcal{F}$ over $K$ and Let $\mathfrak{S}_{n}$ be the symmetric group on $n$ letters. We always denote $(b_{ij})$ for the square matrix $B$ and $[1,n]=\{1,2, \ldots, n \}$.\\
{\textbf{Acknowledgements.} The author is very grateful to the referee for the  comments and suggestions which were very useful  in finishing the paper in the final form.
\section{Cluster algebras associated with valued quivers}
For more details about the material of this section  refer to  [4, 9, 5, 6, 7].
\subsection{Valued quivers}
 \begin{itemize}
\item
   \emph{A valued quiver} of rank $n$ is a quadruple $Q=(Q_{0}, Q_{1}, V, d)$, where
  \begin{itemize}
 \item $Q_{0}$ is a  set of $n$ vertices labeled by $[1,n]$;
  \item  $Q_{1}$ is called the \emph{set of arrows} of $Q$ and consists of  ordered pairs of vertices, that is $Q_{1}\subset Q_{0}\times Q_{0}$;
  \item $V=\{(v_{ij},v_{ji})\in \mathbb{N}\times\mathbb{N} | (i,j)\in Q_{1}\}$, $V$ is called the \emph{valuation} of $Q$;
  \item $d=(d_{1},\cdots, d_{n})$, where $d_{i}$ is a positive integer for each $i$, such that $d_{i}v_{ij}=v_{ji}d_{j}$, for every $i, j\in [1, n]$.
   \end{itemize}
   In the case of $(i,j)\in Q_{1}$, then there is an arrow oriented from $i$ to $j$ and in notation we shall use the symbol $\xymatrix{{\cdot}_{i} \ar[r]^{(v_{ij},v_{ji})}&{\cdot}_{j}}$. If $v_{ij}=v_{ji}=1$ we simply write $\xymatrix{{\cdot}_{i} \ar[r]&{\cdot}_{j}}$.\\
In this paper, we moreover assume that $(i, i)\notin Q_{1}$ for every $i\in Q_{0}$, and  if $(i, j)\in Q_{1}$ then $(j, i)\notin Q_{1}$
 \item   If $v_{ij}=v_{ji}$ for every $(v_{ij},v_{ji})\in V$ then $Q$ is called  \emph{equally  valued quiver}.

  \item A \emph{valued quiver morphism} $\phi$ from  $Q=(Q_{0}, Q_{1}, V, d)$ to $Q'=(Q'_{0}, Q'_{1}, V', d')$ is a pair of maps $(\sigma _{\phi}, \sigma _{1})$ where $\sigma _{\phi}:Q_{0}\rightarrow Q_{0}'$ and $\sigma _{1}:Q_{1}\rightarrow Q_{1}'$ such that  $\sigma_{1}(i, j)=(\sigma_{\phi}(i), \sigma_{\phi}(j)) $ and $(v_{ij}, v_{ji})=(v'_{\sigma _{\phi}(i)\sigma _{\phi}(j)}, v'_{\sigma _{\phi}(j)\sigma _{\phi}(i)})$ for each  $(i, j)\in Q_{1}$. If $\phi$ is invertible then it is called a \emph{valued quiver isomorphism}. In particular $\phi$ is called a \emph{quiver automorphism} of $Q$ if it is a valued quiver isomorphism from $Q$ to itself.
 \end{itemize}
 \begin{rems}
 \begin{enumerate}
   \item
   Every (non valued) quiver $Q$ without loops nor $2$-cycles corresponds to an equally valued quiver which has an arrow $(i, j)$ if there is at least one arrow directed from $i$ to $j$ in $Q$ and with the valuation   $(v_{ij}, v_{ji})=(m, m)$, where $m$ is the number of arrows from $i$ to $j$.

  \item
Every  valued quiver of rank $n$ corresponds to a skew symmetrizable integer matrix   $B(Q)=(b_{ij})_{i,j\in[1,n]}$ given by
\begin{equation}\label{}
  b_{ij}=\begin{cases} v_{ij}, & \text{ if }(i,j)\in Q_{1},\\
    0, & \text{if \  neither} \  (i,j)  \ \text{nor} \  (j,i) \  \text{is in} \  Q_{1},\\
-v_{ij}, & \text{ if }(j,i)\in Q_{1}.
    \end{cases}
\end{equation}
Conversely,  given a skew symmetrizable $n\times n$ matrix $B$, a valued quiver $Q_{B}$ can be easily defined  such that $B(Q_{B})=B$. This gives rise to a  bijection between the skew-symmetrizabke  $n\times n$ integral matrices $B$ and the valued quivers with set of vertices $[1,n]$, up to isomorphism fixing the vertices.

  \end{enumerate}
 \end{rems}

\begin{itemize}
\item
The \emph{mutation of valued quivers }is defined through  Fomin-Zelevinsky's mutation of the associated skew-symmetrizable matrix. The mutation of a skew symmetrizable matrix $B=(b_{ij})$ on the direction $k\in[1,n]$ is given by $\mu_{k}(B)=(b'_{ij})$, where
\begin{equation}\label{}
b'_{ij}=\begin{cases} -b_{ij}, & \text{ if }k\in \{i,j\},\\
    b_{ij}+\text{sign}(b_{ik})\max(0, b_{ik}b_{kj}), & \text{otherwise.}
    \end{cases}
  \end{equation}
  \item
  The mutation of a valued quiver $Q$ can be described using the mutation of $B(Q)$ as follows: Let $\mu _{k}(Q)$ be the valued quiver obtained from  $Q$ by applying mutation at the vertex $k$. We obtain $Q'_{1}$ and $V'$, the set of the arrows  and the valuation of  $\mu _{k}(Q)$ respectively,  by altering  $Q_{1}$ and $V$ of $Q$, based on the following rules
\begin{enumerate}
 \item replace the  pairs $(i, k)$ and $(k,j)$  with $(k,i)$ and $(j,k)$  respectively and switch the components of the ordered pairs of their valuations;
  \item if  $(i,k), (k,j)\in Q_{1}$, but $(j,i)\notin Q_{1}$  and $(i,j)\notin Q_{1}$ (respect to $(i,j)\in Q_{1}$) add the  pair $(i, j)$ to $Q'_{1}$, and give it the valuation $(v_{ik}v_{kj},v_{ki}v_{jk})$ (respect to change its valuation to $(v_{ij}+v_{ik}v_{kj},v_{ji}+v_{ki}v_{jk})$);
 \item if $(i,k)$, $(k,j)$ and $(j, i)$ in $Q_{1}$, then we have three cases
 \begin{enumerate}
   \item if $v_{ik}v_{kj}<v_{ij}$, then keep $(j,i)$ and change its valuation to $(v_{ji}-v_{jk}v_{ki}, -v_{ij}+v_{ik}v_{kj})$;
   \item if $v_{ik}v_{kj}>v_{ij}$, then replace $(j,i)$ with $(i,j)$ and change its valuation to $( -v_{ij}+v_{ik}v_{kj}, |v_{ji}-v_{jk}v_{ki}|)$;
   \item if $v_{ik}v_{kj}=v_{ij}$,  then remove $(j,i)$ and its valuation.
 \end{enumerate}
\end{enumerate}
 \item One can see that; $\mu^{2}_{k}(Q)=Q$ and  $\mu_{k}(B(Q))=B(\mu_{k}(Q))$ for each vertex $k$.
 \end{itemize}

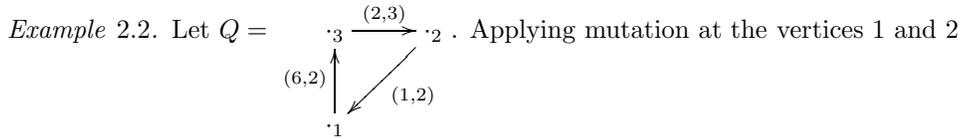
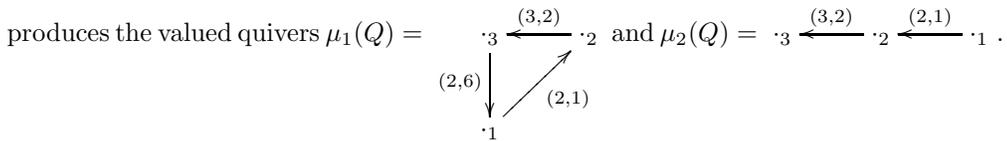
  \begin{exam} Let $Q=\xymatrix{
 \cdot_{3}   \ar[r]^{(2,3)}&   \cdot_{2} \ar[dl]^{(1,2)}\\
  \cdot_{1} \ar[u]^{(6,2)} }$. Applying mutation at the vertices  $1$ and  $2$ produces the valued quivers
  $\mu_{1}(Q)=\xymatrix{ \cdot_{3}   \ar[d]_{(2,6)}&   \cdot_{2}\ar[l]_{(3,2)}\\
  \cdot_{1}\ar[ur]_{(2,1)}  }$ and  $\mu_{2}(Q)=\xymatrix{\cdot_{3}   &   \cdot_{2} \ar[l]_{(3,2)} &\cdot_{1} \ar[l]_{(2,1)} }$.
  \end{exam}

\subsection{Cluster algebras}
 [7, Definition 2.3] \emph{A  labeled seed} of rank $n$ in $\mathcal{F}$  is a pair $(X, Q)$ where $X=(x_{1},\ldots, x_{n})$ is an $n$-tuple elements of $\mathcal{F}$ forming a free generating set and $Q$ is a valued quiver of rank $n$. In this case, $X$ is called a \emph{cluster}.\\ We will  refer to labeled seeds simply as seeds, when there is no risk of confusion.\\
 The  definition of  clusters above  is a bit different from the definition of clusters given in [3], [5] and [6].

\begin{defn}[Seed mutations] Let $p=(X, Q)$ be a seed in $\mathcal{F}$, and $k\in [1,n]$. A new seed $\mu_{k}(X,Q)=(\mu _{k}(X), \mu _{k}(Q))$ is obtained from $(X, Q)$ by  setting $\mu _{k}(X)=(x_{1}, \ldots,x'_{k},\ldots, x_{n})$  where $x'_{k}$ is defined by the so-called \emph{exchange relation}:

\begin{equation}\label{}
   x'_{k}x_{k}=f_{p, x_{k}},\ \ \text{where} \ \ f_{p, x_{k}}=\prod_{(i,k)\in Q_{1}} x_{i}^{v_{ik}}+\prod_{(k,i)\in Q_{1}} x_{i}^{v_{ki}}.
\end{equation}
And $\mu_{k}(Q)$ is the mutation of $Q$ at the vertex $k$. The elements of $\mathcal{F}$ obtained by applying iterated mutations on elements of $X$ are called \emph{cluster variables}.
   \end{defn}

\begin{defns}[Mutation class and cluster algebra]
\begin{itemize}
  \item  The  equivalence class of a seed $(X, Q)$ under mutation is called the \emph{mutation class} of $(X,Q)$ and it will be denoted by Mut$(X, Q)$.
  \item Let $\mathcal{X}$ be the  union of all clusters in Mut$(X, Q)$. The (coefficient free) \emph{cluster algebra} $\mathcal{A}=\mathcal{A}(X, Q)$ is  the $\mathbb{Z}$-subalgebra of $\mathcal{F}$ generated by $\mathcal{X}$.
\end{itemize}
\end{defns}

\begin{thm} [5, Theorem 3.1, Laurent Phenomenon]  The cluster algebra
$\mathcal{A}(X, Q)$ is contained in the integral ring of Laurent polynomials $\mathbb{Z}[x^{\pm}_{1},\ldots ,x^{\pm}_{n}]$. More precisely, every non zero element $y$ in $\mathcal{A}(X, Q)$
 can be uniquely written as
\begin{equation}\label{}
    y=\frac{P(x_{1}, x_{2},\ldots ,x_{n})}{x_{1}^{\alpha _{1}}\cdots x_{n}^{\alpha
    _{n}}},
\end{equation}
where $\alpha_{1}, \ldots,\alpha_{n}$ are integers  and $P(x_{1}, x_{2},\ldots, x_{n})$ is a  polynomial with integer coefficients   which is not
divisible by any of the cluster variables $x_{1},\ldots, x_{n}$.
 \end{thm}

\begin{conj}[\emph{Positivity Conjecture}] If $y$ is a cluster variable, then the polynomial $P(x_{1}, x_{2},\ldots, x_{n})$, in (2.4), has positive integer coefficients.
\end{conj}
\begin{defn} \emph{Positive cluster algebras} are the cluster algebras that   satisfy the positivity conjecture.
\end{defn}

\section{isomorphisms of cluster algebras}

\subsection{Quiver automorphisms  action on seeds}
 Let $\phi$ be a quiver automorphism of  $Q =(Q_{0}, Q_{1}, V, d)$. Then $\phi$ induces a permutation $\sigma _{\phi}\in \mathfrak{S} _{n}$. We can obtain a new valued quiver  $\phi (Q)=(Q_{0}', Q_{1}', V', d')$  from $Q$ as follows
 \begin{itemize}
   \item $Q_{0}'$ is obtained by permuting the vertices of $Q_{0}$  using $\sigma_{\phi}$;
   \item $Q_{1}'=\{(\sigma_{\phi} (i), \sigma _{\phi}(j))| (i,j)\in Q_{1} \}$;
   \item For every $(\sigma_{\phi} (i), \sigma_{\phi} (j))\in Q'_{1}$ we give the valuation $(v_{\sigma_{\phi} (i)\sigma_{\phi} (j)},v_{\sigma_{\phi} (j)\sigma_{\phi} (i)})$;
   \item $d'=(d_{\sigma_{\phi}(1)},\cdots, d_{\sigma_{\phi} (n)})$.
 \end{itemize}

\begin{exam} Consider the valued quiver
$Q=\xymatrix{
 \cdot_{1} \ar[d]_{(4,1)}  \ar[r]^{(2,1)}&\cdot_{2} \\
  \cdot_{3}  \ar[ur]_{(1,2)}}$ with $d=(1,2,4)$ and the quiver automorphism $\phi$ with underlying permutation $\sigma_{\phi}=(123)$. Then $\phi(Q)=\xymatrix{
 \cdot_{2} \ar[d]_{(4,1)}  \ar[r]^{(2,1)}&\cdot_{3} \\
  \cdot_{1}  \ar[ur]_{(1,2)}}$ and $\phi (d)=(4, 1, 2)$.
\end{exam}

\begin{que}
For which  valued quivers $Q$  and a quiver automorphism  $\phi$ are  there  a sequence of mutations $\mu$ such that $\phi (Q)=\mu (Q)$?

\end{que}
 The following example shows that there are cases in which such a sequence of mutations does not exist.

\begin{exam} Consider the  valued quiver $Q=\xymatrix{\cdot_{1} \ar[r]^{(2,2)}  &   \cdot_{2} \ar[r] &\cdot_{3}  }$, and the quiver automorphism with underlying permutation $(12)$. In the following, we will show that there is no sequence of mutations
$\mu$, such that $\mu(Q)=\xymatrix{\cdot_{2} \ar[r]^{(2,2)}  &   \cdot_{1} \ar[r] &\cdot_{3}  }$ or equivalently there is no sequence of mutations $\mu $ such that
$\mu(B(Q))=\left( \begin{array}{ccc}
    0 & -2 & +1 \\
    +2 & 0 & 0 \\
    -1 & 0 & 0 \\
  \end{array}
\right)$. Here we have   $B(Q)=  \left( \begin{array}{ccc}
    0 & +2 & 0 \\
    -2 & 0 & +1 \\
    0 & -1 & 0 \\
  \end{array}
\right)$.
The proof is written in terms of the matrix $B(Q)=(b_{ij})$. If we could show that there is no sequence of mutations that
sends the entry $b_{23}$ to zero, we will be done. We do this by showing  that every sequence of
mutations sends $b_{23}$ to an odd number.  First we show, by induction on the length of the sequence of mutations, that any sequence of mutations sends $b_{13}$
and $b_{12}$ to  even numbers. For sequences containing only one mutation: one can see that, only $\mu _{2}$
and $\mu _{3}$ would change $b_{13}$ and $b_{12}$ respectively, that is $\mu _{2}(b_{13})=\mu _{3}(b_{12})=2$.\\
Now, assume that every sequence of mutations of length $k$ sends
$b_{13}$  and $b_{12}$ to an even number, and let $\mu
_{i_{k+1}}\mu _{i_{k}}\ldots \mu _{i_{1}}$ be a sequence of length
$k+1$. So if

\begin{equation}\label{}
    \mu _{i_{k}}\ldots \mu _{i_{1}}((b_{ij}))=(b'_{ij}).
\end{equation}
 Then
$b'_{12}=2m$ for some integer number $m$. We have

\begin{eqnarray}
  \nonumber \mu_{i_{k+1}}(b'_{13})&=& b'_{13}+\text{sign}(b'_{12})\max(0, b'_{12}b'_{23}) \\
  \nonumber &=& b'_{13}+\text{sign}(b'_{12})\max(0, 2mb'_{23})
  \end{eqnarray}
   which is a sum of two even numbers. This shows that any sequence of mutations
will send $b_{13}$ to an even number. In a  similar way one
can show that any sequence of mutation sends $b_{12}$ to an even
number.\\
Secondly, we show that every sequence of mutations sends
$|b_{23}|$ to an odd number, noting that the possible change in $|b_{23}|$ appears only after applying $\mu _{1}$. We show this by induction on the
number of occurrences of $\mu _{1}$ in the sequence.\\
\textbf{Sequences containing only one copy of $\mu _{1}$:} Without loss of generality, let $\mu_{i_{1}}\mu_{i_{2}}\ldots \mu_{i_{k}}$ be a sequence of mutations such that $\mu _{i_{k}}=\mu _{1}$, and $\mu
_{i_{j}}\neq \mu _{1}, \ \ \forall j \in [1, k-1]$. Then using the same notation as in (3.1), we have
\begin{equation}\label{}
    b'_{23}=\pm 1+\text{sign}(b'_{21})\max(0, b'_{21}b'_{13}).
\end{equation}
However  $b'_{21}$  and $b'_{13}$ are both even numbers, so
$\text{sign}(b'_{21})\max(0, b'_{21}b'_{13})$ must be an even number, and then $b'_{23}$ is an odd number.\\
 \textbf{Sequences containing more than one copy of $\mu _{1}$:}
Assume that any sequence of mutations with $\mu _{1}$ repeated
$k$ times sends $b_{23}$ to an odd number.\\
Let $\mu_{i_{t}}\mu_{i_{2}}\ldots \mu_{i_{1}}$ be a sequence of
mutations containing $k+1$-copies of $\mu _{1}$. Then we can assume
that  $\mu_{i_{t}}=\mu_{1}$. Let
\begin{equation}\label{}
    \mu _{i_{t}}\ldots \mu _{i_{1}}((b_{ij}))=(b''_{ij})  \  \text{and} \ \  \mu _{i_{t-1}}\ldots \mu
    _{i_{1}}((b_{ij}))=(b'_{ij}),
\end{equation}
then one can see that $b'_{23}$ is an odd number and
$b'_{12}$ and $b'_{13}$ are both even numbers. Then,

\begin{equation}\label{}
b''_{23}=b'_{23}+\text{sign}(b'_{13})\max(0, b'_{21}b'_{13}),
\end{equation}
a sum of an odd and an even number, hence $b''_{23}$ is
an odd number.                  \ \ \ \ \ \ \ \ \ \ \ \ \ \ \ \ \ \ \ \ \ \ \ \ \ \ \ \          $\Box$
\end{exam}

Let $T=\{t_{1}, t_{2},\ldots, t_{n}\}$ be a free generating set of $\mathcal{F}$ and $\sigma \in \mathfrak{S}_{n}$. We  define an $\mathcal{F}$-automorphism $\sigma _{T}$, by  $\sigma _{T}(r(t_{1},\ldots t_{n}))=r(t_{\sigma(1)}, \ldots t_{\sigma (n)})$ for  $r(t_{1}, \ldots t_{n})\in \mathcal{F}$.

\begin{defn} Fix a cluster $X$ in $\mathcal{F}$. A quiver automorphism $\phi$, with underlying permutation  $\sigma_{\phi} \in \mathfrak{S}_{n}$, acts on the set of all seeds of $\mathcal{F}$  \emph{with respect to} $X$  as follows:
 for any seed  $q=(Y, \Gamma)$ define
\begin{equation}\label{}
    \phi _{X}(q):=(\phi_{X}(Y), \phi(\Gamma)),
\end{equation}
where $\phi _{X}(Y)=((\sigma_{\phi})_{ X}(y_{1}), \ldots, (\sigma _{\phi})_{ X}(y_{n}))$. We write   $\phi(p)$ and $\phi (Y)$ instead of $\phi _{X}(q)$ and $\phi_{X}(Y)$ respectively if there is no chance of confusion.
\end{defn}

\begin{defns}\begin{enumerate}
\item
 Let  $\sigma \in \mathfrak{S}_{n}$. Two valued quivers $Q$ and $Q'$ are said to be \emph{$\sigma$-similar} if $\sigma $ is the underlying permutation of  a quiver isomorphism between $Q$ and one of the valued quivers $Q'$ or $(Q')^{\text{op}}$.
 \item  Two seeds $(X, Q)$ and $(Y, Q')$ are said to be $\sigma$\emph{-similar} if $Q$ and $Q'$ are $\sigma$-similar.
\end{enumerate}
\end{defns}

\begin{rem}\begin{enumerate}
             \item [(a)]$Q$ and $Q'$ are $\sigma$-similar if and only if $B(Q)=\epsilon \sigma ( B(Q'))$, for $\epsilon \in \{-1, +1\}$.
             \item  [(b)] The $\sigma$-similarity relation defines an equivalence relation on the set of all seeds of $\mathcal{F}$.
           \end{enumerate}

\end{rem}

\begin{lem}
Two quivers $Q$ and $Q'$ are $\sigma$-similar if and only if
\begin{equation}\label{}
\prod_{(i,k)\in Q_{1}} t_{i}^{v_{ik}}+\prod_{(k,i)\in Q_{1}} t_{i}^{v_{ki}}=\prod_{(\sigma (i),\sigma (k))\in Q_{1}} t_{\sigma (i)}^{v_{\sigma (i)\sigma (k)}}+\prod_{(\sigma (k),\sigma (i))\in Q_{1}} t_{\sigma (i)}^{v_{\sigma (k)\sigma (i)}},  \forall k\in [1, n]
\end{equation}

\end{lem}
\begin{proof} One can see that $Q$ and $Q'$ are $\sigma$-similar if and only if  $(v_{ij})_{i,j \in [1, n]}=(v'_{\sigma(i)\sigma (j)})_{i,j \in [1, n]}$  and one of the following two conditions is satisfied

  \begin{equation}\label{}
     (i, j)\in Q_{1} \ \ \text{if and only if} \ \  (\sigma (i), \sigma (j))\in Q'_{1}, \   \text{for every} \ \ i, j\in [1, n];
  \end{equation}
  or
 \begin{equation}\label{}
     (i, j)\in Q_{1}\ \ \text{ if and only if } (\sigma (j), \sigma (i))\in Q'_{1}, \ \ \text{for every } \ \ i, j\in [1, n].
 \end{equation}
Since $\{t_{1},\ldots,t_{n}\}$ is a transcendance basis of $\mathcal{F}$, then one of the  conditions  (3.7) or (3.8) is satisfied if and only if one of the following two conditions is satisfied
\begin{enumerate}
  \item

  \begin{equation}\label{}
 \nonumber \prod_{(i,k)\in Q_{1}} t_{i}^{v_{ik}}=\prod_{(\sigma(i),\sigma(k))\in Q_{1}} t_{\sigma (i)}^{v_{\sigma(i)\sigma(k)}}\  \text{and} \ \prod_{(k,i)\in Q_{1}} t_{i}^{v_{ki}}=\prod_{(\sigma(k),\sigma (i))\in Q_{1}} \ t_{i}^{v_{\sigma(k)\sigma (i)}},   \forall k\in [1, n];
\end{equation}
 or
\item

 \begin{equation}\label{}
 \nonumber \prod_{(i,k)\in Q_{1}} t_{i}^{v_{ik}}=\prod_{(\sigma(k),\sigma (i))\in Q_{1}} \ t_{i}^{v_{\sigma(k)\sigma (i)}}\  \text{and} \ \prod_{(k,i)\in Q_{1}} t_{i}^{v_{ki}}=\prod_{(\sigma(i),\sigma(k))\in Q_{1}}t_{\sigma (i)}^{v_{\sigma(i)\sigma(k)}},   \forall k\in [1, n].
\end{equation}
\end{enumerate}
Which is equivalent to (3.6).

\end{proof}

\subsection{Main definitions}
 Let $p=(X, Q)$ and $p'=(Y, Q')$ be two seeds of rank $n$.
    \begin{defn}
Let $f$ be an element of $Aut_{K}(\mathcal{F})$.
\begin{itemize}
\item  $f$ is called a \emph{cluster variables preserver} from $\mathcal{A}(X,Q)$ to $\mathcal{B}(Y,Q')$, if  it sends every cluster variable in $\mathcal{A}$ to a cluster variable in $\mathcal{B}$.
   In particular, $f$ is  a cluster variables preserver of a cluster algebra $\mathcal{A}(X, Q)$ if it leaves $\mathcal{X}$, the set of all cluster variables of $\mathcal{A}$, invariant. For simplicity we  call such automorphisms
    the $\mathcal{X}$-\emph{preservers}.

 \item \([1, 2, 6]\) $f$ is said to be a \emph{cluster isomorphism} from $\mathcal{A}(X, Q)$ to $\mathcal{B}(Y, Q')$ if it induces a
one to one correspondence from the set of all clusters of $\mathcal{A}$ to the set of all clusters of $\mathcal{B}$. In particular, $f:\mathcal{A} \rightarrow \mathcal{A}$
is called a \emph{cluster automorphism} of $\mathcal{A}$ if it permutes the clusters of $\mathcal{A}$.
\end{itemize}
\end{defn}
\begin{rem}
 An element  $f$ of $Aut_{K}(\mathcal{F})$ is a cluster isomorphism from  $\mathcal{A}(X, Q)$ to $\mathcal{B}(Y, Q')$ if and only if it induces a one to one correspondence
between the mutation classes  Mut$(X, Q)$ and  Mut$(Y, Q')$. This  is due to the fact that  for every seed $(X, Q)$ the quiver $Q$ is uniquely defined by the cluster $X$ which has been proved in [8, Theorem 3].

\end{rem}

\subsection{Exchange maps}

\begin{defn} Let $\sigma\in \mathfrak{S}_{n}$ and $p=(X, Q)$ and $p'=(Y, Q')$ be two labeled seeds in $\mathcal{F}$. Then, the map   $T_{pp', \sigma}$, induced by $x_{i}\mapsto y_{\sigma (i)}$, is called an \emph{exchange map}.
\end{defn}
One can see that exchange maps are elements of $Aut_{K}(\mathcal{F})$.
\begin{lem} Let  $\sigma \in \mathfrak{S}_{n}$. Then
\begin{equation}\label{}
  \nonumber Q  \ \ \text {and} \ \  Q' \ \ \text {are} \ \ \sigma-\text{similar} \ \ \text{if and only
    if} \ \  T_{pp', \sigma}(\mu _{k}(x_{k}))=\mu _{\sigma (k)}(y_{\sigma
    (k)}),\ \ \text{for all} \  k \in [1, n].
\end{equation}
\begin{proof}$\Rightarrow$) Assume that $Q$ and $Q'$ are  $\sigma$- similar valued quivers. Then  $v_{ij}=v'_{\sigma (i)\sigma (j)}$ for every $i, j \in [1, n]$ and one of the conditions (3.7) or (3.8)  is satisfied. Hence
\begin{eqnarray}
 \nonumber T_{pp', \sigma} ( \mu _{k}(x_{k}))&=&T_{pp', \sigma} \left(\frac{1}{x_{k}}\left(\prod_{(i, k)\in Q_{1}}x^{v_{ik}}_{i}+\prod_{(k, i)\in Q_{1}}x^{v_{ki}}_{i}\right)\right)\\
   \nonumber&=& \frac{1}{y_{\sigma (k)}}\left(\prod_{(i, k)\in Q_{1}}y^{v_{ik}}_{\sigma (i)}+\prod_{(k, i)\in Q_{1}}y^{v_{ki}}_{\sigma (i)}\right)\\
   \nonumber&=& \frac{1}{y_{\sigma (k)}}\left(\prod_{(\sigma (i), \sigma (k))\in Q_{1}'}y^{v'_{\sigma (i) \sigma (k)}}_{\sigma (i)}+\prod_{(\sigma (k), \sigma (i))\in Q_{1}'}y^{v'_{\sigma (k) \sigma (i)}}_{\sigma (i)}\right)\\
   \nonumber &=&\mu _{\sigma (k)}(y_{\sigma (k)}).
\end{eqnarray}
$\Leftarrow$) Suppose that $Q$ and $Q'$ are not $\sigma$-similar. Then (3.6) is not satisfied. Hence $T_{pp', \sigma}(f_{p, x_{k}})\neq f_{p', y_{\sigma (k)}}$, for some  $k\in [1, n]$. Therefore $T_{pp', \sigma}(\mu _{i}(x_{k}))\neq \mu _{\sigma (k)}(y_{\sigma (k)})$ for some $k\in [1, n]$.
\end{proof}

\end{lem}
\begin{lem} For every $\sigma \in \mathfrak{S_{n}}$ and any square matrix $B$, we have
\begin{equation}\label{}
 \sigma (\mu_{k}(B))=\mu_{\sigma (k)}(\sigma(B)), \ \ \text{for all} \ k \in  [1, n].
\end{equation}
In particular, for every valued quiver $Q$
\begin{equation}\label{}
     \sigma (\mu_{k}(Q))=\mu_{\sigma (k)}(\sigma(Q)), \ \ \ \text{for all} \in  [1, n].
\end{equation}
\begin{proof} Let $\sigma (B)=(b^{*}_{ij})$, $\mu_{\sigma (k)} (\sigma (B))=(\overline{b}_{ij})$, $\mu_{k}(B)=(b'_{ij})$, and  $\sigma (\mu_{k}(B))=(b^{\star}_{ij})$. We obtain the matrix $\sigma (B)$ from $B$,  by relocating the entries of $B$ using $\sigma$. Indeed, the entry $b^{*}_{ij}=b_{\sigma^{-1} (i)\sigma^{-1} (j)}$.
\begin{eqnarray*}
  \nonumber \overline{b}_{ij}&=&\begin{cases} -b^{*}_{ij}, & \text{ if }  \sigma(k) \in \{i,j\},\\
       b^{ *}_{ij}+\text{sign}(b^{*}_{ik})\text{max}(0, b^{*}_{ik}b^{*}_{kj}), & \text{otherwise}
    \end{cases}    \\
  \nonumber &=&\begin{cases}
     -b_{\sigma^{-1}(i)\sigma^{-1}(j)},  \ \ \ \ \ \ \ \ \ \ \  \ \ \ \ \ \ \  \ \ \ \text{ if } k \in \{\sigma^{-1}(i),\sigma^{-1}(j)\}&\\
         b_{\sigma^{-1}(i)\sigma^{-1}(j)}+\text{sign}(b_{\sigma^{-1}(i)\sigma^{-1}(k)})\text{max}(0, b_{\sigma^{-1}(i)\sigma^{-1}(k)}b_{\sigma^{-1}(k)\sigma^{-1}(j)}), & \text{otherwise}
    \end{cases}    \\
   \nonumber &=& b'_{\sigma ^{-1}(i)\sigma ^{-1}(j)}\\
 \nonumber &=& b^{\star}_{ij}.\\
\end{eqnarray*}
This proves (3.9). Identity (3.10)  is immediate by using $B=B(Q)$ in (3.9).

\end{proof}
\end{lem}

\begin{thm} Let $(X, Q)$ and $(Y, Q')$ be two $\sigma$-similar seeds. Then for any sequence of mutations $\mu _{i_{k}}, \mu _{i_{k-1}}, \ldots, \mu _{i_{1}}$, the following are true
\begin{enumerate}

\item
 $\mu _{i_{k}}\mu _{i_{k-1}}\ldots\mu _{i_{1}}(X, Q)$ and $\mu _{\sigma (i_{k})}\mu _{\sigma (i_{k-1})}\ldots\mu _{\sigma (i_{1})}(Y, Q')$ are $\sigma$- similar,
\item
 $ T _{pp',\sigma}(\mu _{i_{k}}\mu _{i_{k-1}}\ldots\mu _{i_{1}}(X))=\mu _{\sigma (i_{k})}\mu _{\sigma (i_{k-1})}\ldots\mu _{\sigma (i_{1})}(Y)$.
\end{enumerate}
\begin{proof}
To prove part $(1)$, assume that $Q$ and $Q'$ are two $\sigma$-similar valued quivers. Then $B(Q)=\epsilon \sigma ( B(Q'))$ for $\epsilon \in \{-1, +1\}$. Hence from   (3.9), we have

\begin{equation}\label{}
 \nonumber  \mu_{k}(B(Q))= \mu_{k}(\epsilon \sigma (B(Q'))) = \epsilon \sigma (\mu_{\sigma (k)}( B(Q'))), \ \text{for every}\ k \in[1, n].
\end{equation}
 Therefore, $\mu_{k}(B(Q))$ and   $\mu_{\sigma (k)} (B(Q'))$ are $\sigma$-similar  for every $k\in [1, n]$. An induction process generalizes this fact to any  sequence of mutations $\mu _{i_{k}}, \mu _{i_{k-1}}, \ldots, \mu _{i_{1}}$.\\
 For part $(2)$, let $p_{i_{1}i_{k}}=\mu _{i_{k}}\mu _{i_{k-1}}\ldots\mu _{i_{1}}(p)$ and $p'_{\sigma (i_{1})\sigma (i_{k})}=\mu _{\sigma (i_{k})}\mu _{\sigma (i_{k-1})}\ldots\mu
 _{\sigma (i_{1})}(p')$. Part (1) of this theorem tells us that $p_{i_{1}i_{k-1}}$ and  $p'_{\sigma (i_{1})\sigma (i_{k-1})}$ are $\sigma$-similar.
 Then Lemma 3.11 implies that

 \begin{equation}\label{}
   T_{p_{i_{1}i_{k}}p'_{\sigma (i_{1})\sigma (i_{k})}, \sigma } (\mu _{i_{k}}(\mu _{i_{k-1}}\ldots\mu
    _{i_{1}}(X)))= \mu _{i_{\sigma (k)}}(\mu _{i_{\sigma (k-1)}}\ldots\mu _{\sigma (i_{1})}(Y)).
 \end{equation}
So, it remains to show that

 \begin{equation}\label{}
      T_{p_{i_{1}i_{k}}p'_{\sigma (i_{1})\sigma (i_{k})}, \sigma } (\mu _{i_{k}}(\mu _{i_{k-1}}\ldots\mu
       _{i_{1}}(X)))= T_{pp', \sigma } (\mu _{i_{k}}(\mu _{i_{k-1}}\ldots\mu_{i_{1}}(X)).
  \end{equation}
 This will be proved by induction on the length of the sequence of mutations. First we will show (3.12) in the general setting for sequences of mutations of length 2. Let $q=(Z, D)$ and $q'=(T, C)$ be any two $\sigma$-similar  seeds, and let $q_{1}=\mu _{i}(Z, D)=(Z', D'), \ \ q'_{1}= \mu_{\sigma (i)}(T, C)=(T', C')$. We  show that

\begin{equation}\label{}
T_{q_{1}q'_{1}, \sigma}(\mu_{k}\mu_{i}(Z))=T_{qq', \sigma}(\mu_{k}\mu_{i}(Z)).
\end{equation}
Let $z_{j}$ be a cluster variable in $Z$. Then for $j\neq i$, both of
$T_{q_{1}q'_{1}, \sigma}$, and $T_{qq', \sigma}$ will leave $z_{j}$ unchanged. Now, let $j=i$. Then, using the notation $f_{p, x_{k}}$ introduced in Equation (2.3) above, we have

\begin{equation}\label{}
    \nonumber T_{q_{1}q'_{1}, \sigma}(\mu _{i}(z_{i}))=T_{q_{1}q'_{1}, \sigma} \left(\frac{f_{q,z_{i}}}{z_{i}}\right ) = \frac{f_{q'_{1},t_{\sigma (i)}}}{T_{p_{1}p'_{1}, \sigma}(z_{i})}.
 \end{equation}
However,

\begin{equation}\label{}
   \nonumber T_{q_{1}q'_{1}, \sigma}(\mu _{i}(z_{i}))= \mu _{\sigma (i)}(t_{\sigma(i)})=  \frac{f_{q'_{1},t_{\sigma (i)}}}{t_{\sigma(i)}}.
    \end{equation}
Hence, $T_{p_{1}p'_{1}, \sigma}(z_{i})=t_{\sigma (i)}$. This shows that   $ T_{q_{1}q'_{1}, \sigma}$, and $T_{qq', \sigma}$ have the same action
 on every cluster variable in $Z$, and since  cluster variables from the cluster $\mu _{k}\mu _{i}(Z)$ are  integral Laurent polynomials
 of cluster variables from $Z$, this gives (3.13).\\
For equation (3.12) we use induction on the length of the mutation sequence. Assume that
equation (3.12) is true for any sequence of mutations of length less than or equal $k-1$. Now we have;

\begin{eqnarray}
  \nonumber T _{p_{i_{1}i_{k}}p'_{i_{1}i_{k}}, \sigma }(\mu _{i_{k}}\mu _{i_{k-1}} \ldots \mu _{i_{1}}(X))&=& T _{p_{i_{1}i_{k-2}}p'_{i_{1}i_{k-2}}, \sigma }(\mu _{i_{k}}\mu _{i_{k-1}} \ldots \mu _{i_{1}}(X)) \\
   \nonumber&=& T_{pp',\sigma}(\mu _{i_{k}}\mu _{i_{k-1}} \ldots \mu _{i_{1}}(X)),
  \nonumber\end{eqnarray}
where the first equality is by  (3.13) and the second is  by the induction hypotheses.

\end{proof}

\end{thm}

\subsection{Main theorem}

Equivalent conditions for exchange maps to be cluster automorphisms are provided in this subsection.

\begin{thm} Let $p=(X, Q)$ and $p'=(Y, Q')$ be two labeled seeds in $\mathcal{F}$, and $\sigma\in \mathfrak{S}_{n}$. Then $p$ and $p'$ are $\sigma$-similar if and only if $T _{pp', \sigma}$ is a cluster isomorphism from $\mathcal{A}(X, Q)$ to $\mathcal{B}(Y, Q')$. In particular, two  cluster algebras of the same rank $n$ are cluster isomorphic if and only if there exists a permutation $\sigma\in \mathfrak{S}_{n}$ such that the cluster algebras  contain two $\sigma$-similar seeds.
\end{thm}
\begin{proof} Assume that $p$ and $p'$ are $\sigma$-similar. Let $(Z,D) \in$ Mut$(X, Q)$. Then there is a sequence of mutations $\mu _{i_{1}},\mu
_{i_{2}},\ldots, \mu _{i_{k}} $ such that $Z=\mu _{i_{1}}\mu _{i_{2}}\ldots \mu_{i_{k}}(X)$. Then part $(1)$ of  Theorem 3.13  implies that  $\mu _{i_{2}}\ldots \mu _{i_{k}}(X, B)$ and $\mu _{\sigma (i_{2})}\ldots \mu _{\sigma (i_{k})}(Y, B')$ are $\sigma$-similar too. But, from  Theorem 3.13 part $(2)$, we have
\begin{eqnarray}
   \nonumber T_{pp', \sigma}(\mu _{i_{1}}(\mu _{i_{2}}\ldots \mu
   _{i_{k}}(X)))=\mu _{\sigma (i_{1})}(\mu _{\sigma (i_{2})}\ldots \mu _{\sigma
   (i_{k})}(Y)),
\end{eqnarray}
and since the right hand side is a cluster, then $T_{pp', \sigma}$ sends $Z$ to a
cluster in Mut$(Y, Q')$. So, $T_{pp', \sigma}$ sends every cluster in Mut$(X, Q)$ to a cluster in Mut$(Y, Q')$. Since the map $\mu _{i_{1}}\ldots \mu _{i_{t}}\mapsto \mu _{\sigma (i_{1})}\ldots \mu _{\sigma(i_{t})}$, with $t$ a non-negative integer, is a one to one correspondence on the set of all sequences of mutations. Thus $T_{pp', \sigma}$ defines  a one to one correspondence from the set of all clusters of  Mut$(X, Q)$ to the set of all clusters of Mut$(Y, Q')$. \\Assume that $T_{pp', \sigma}$ is a cluster isomorphism.  Then $T_{pp', \sigma}(\mu_{i}(X))$ is a cluster in Mut$(Y, Q')$; which shares  $n-1$ cluster variables with the cluster  $\mu _{\sigma (i)}(Y)$. Then from [8, Theorem 3], each one of the two clusters can be obtained from the other by applying one mutation which must be $\mu_{\sigma(i)}$.  But $\mu^{2}_{\sigma(i)}(y_{\sigma(i)})=y_{\sigma(i)}$  and $T_{pp', \sigma}(\mu_{i}(x_{i}))=\frac{T_{pp', \sigma}(f_{p, x_{i}})}{y_{\sigma(i)}}$. Then the two clusters coincide, hence $T_{pp', \sigma}(\mu _{i}(x_{i}))=\mu _{\sigma (i)}(y_{\sigma (i)})$ for every $i \in [1, n]$. Therefore, $p$ and $p'$ are $\sigma$-similar, thanks to  Lemma 3.11.
\end{proof}
\begin{thm} Let $p=(X, Q)$ and $p'=(Y, Q')$ be two labeled seeds, of the same rank $n$, in $\mathcal{F}$ and $\sigma\in \mathfrak{S}_{n}$. Then, if both of $\mathcal{A}(X, Q)$ and $\mathcal{B}(Y, Q')$ are positive cluster algebras, then the following are
equivalent:

\begin{enumerate}
  \item  $T _{pp', \sigma}$ is a cluster variables preserver from  $\mathcal{A}(X, Q)$ to  $\mathcal{B}(Y, Q')$;
  \item $p$ and $p'$ are $\sigma$-similar;
  \item $T _{pp', \sigma}$ is a cluster isomorphism from $\mathcal{A}(X, Q)$ to $\mathcal{B}(Y, Q')$.
  \end{enumerate}

\begin{proof}
$(1)\Rightarrow (2)$. To show  that $p$ and $p'$ are $\sigma$-similar,  we only need to show that  $T_{pp', \sigma}(\mu _{i}(x_{i}))=\mu _{\sigma (i)}(y_{\sigma (i)}), \ \text{for all} \  i\in [1, n]$, thanks to  Lemma 3.11.\\
 Let $z=T_{pp',\sigma}(\mu _{i}(x_{i}))$ and $\xi=\mu _{\sigma (i)}(y_{\sigma
(i)})$. Then

\begin{equation}\label{}
    z=\frac{T_{pp', \sigma }(f_{p, x})}{y_{\sigma (i)}}, \ \ \  \ \ \  \text
{and} \ \ \  \ \ \  \xi=\frac{f_{p', y_{\sigma (i)}}}{y_{\sigma (i)}}.
\end{equation}
 Both of $T_{pp', \sigma }(f_{p, x})$ and $f_{p', y_{\sigma(i)}}$ are  elements in  the ring of polynomials\\ $\mathbb {Z}[y_{\sigma (1)},\cdots,y_{\sigma
(i-1)},y_{\sigma (i+1)},\ldots, y_{\sigma (n)}]$ such that neither of them is divisible by $y_{\sigma (j)}$ for any $j$ $\in[1, n]$. From the assumption that $T _{pp', \sigma}$ is a cluster variables preserver from cluster algebra $\mathcal{A}(X, Q)$ to the cluster algebra $\mathcal{B}(Y, Q')$, then  $z$ must be a cluster variable in $\mathcal{B}(Y, Q')$. Hence, by Laurent phenomenon,  $z$ is an element of the ring of Laurent polynomials in the variables of the cluster $\mu_{\sigma (i)}(Y)$ with integer coefficients. More precisely, $z$  can be written uniquely as
\begin{equation}\label{}
   z=\frac{P(y_{\sigma (1)}, y_{\sigma (2)},\ldots,y_{\sigma (i-1)}, \xi, y_{\sigma (i+1)}, \ldots, y_{\sigma (n)})}{y_{\sigma (1)}^{\alpha _{1}}\ldots y^{\alpha _{i-1}}_{\sigma (i-1)} \xi^{\alpha _{i}} y^{\alpha _{i+1}}_{\sigma (i+1)}\ldots y_{\sigma (n)}^{\alpha
    _{n}}},
\end{equation}
where $\alpha _{i}\in \mathbb{Z}$, for all $i\in [1, n]$ and $P$ is  not divisible by any of the cluster variables $y_{\sigma (1)}, y_{\sigma(2)}, \ldots y_{\sigma (i-1)}, \xi, y_{\sigma (i+1)}, \ldots, y_{\sigma (n)}$. Comparing $z$ from (3.14) and (3.15),
we get

\begin{equation}\label{}
 T_{pp', \sigma }(f_{p, x})\cdot y_{\sigma (1)}^{\alpha _{1}}\ldots y^{\alpha _{i-1}}_{\sigma (i-1)}\ldots \xi^{\alpha _{i}} \ldots y^{\alpha _{i+1}}_{\sigma (i+1)}\ldots y_{\sigma (n)}^{\alpha _{n}}=P\cdot y_{\sigma(i)}.
\end{equation}
Since $ f_{p, x}$ is not divisible by any cluster variable $x_{i}$ for any $i\in [1, n]$ as well then $T_{pp', \sigma }(f_{p, x})$ is not divisible by  $y_{i} \ \text{for all} \ i \in [1, n]$. More precisely  $T_{pp', \sigma }(f_{p, x})$ is a sum of two monomials in  cluster variables from  $Y\setminus \{y_{\sigma (i)}\}$, with positive exponents. Therefore, $\alpha _{j}=0$ for all $j\in  [1, n]\setminus \{i\}$. Then equation (3.16) is reduced to

\begin{equation}\label{}
 \nonumber T_{pp', \sigma }(f_{p, x}) \left(\frac{f_{p', y_{\sigma (i)}}}{y_{\sigma (i)}}\right)^{\alpha _{i}}=P\cdot y_{\sigma(i)}.
\end{equation}

 For $\alpha _{i}$, we break it down into three cases; (1) if $\alpha_{i} \geq 0$, then $y^{\alpha _{i}+1}_{\sigma(i)}$ divides either $T_{pp', \sigma }(f_{p, x})$ or $f_{p', y_{\sigma (i)}}$ which is a contradiction,  (2) if $\alpha _{i} < -1$, then $y^{-\alpha _{i}-1}_{\sigma(i)}$ divides either $P$ or $f_{p', y_{\sigma (i)}}$ which again is a contradiction, (3) assume $\alpha _{i}=-1$. Hence (3.16) ends up to

\begin{equation}\label{}
   T_{pp', \sigma }(f_{p, x})=P\cdot f_{p', y_{\sigma (i)}},
\end{equation}
where $f_{p', y_{\sigma (i)}}$ is a sum of two monomials in cluster variables from  $Y\setminus \{y_{\sigma (i)}\}$, with positive exponents. However, $P$ is a polynomial with positive integers coefficients, and not divisible by any cluster variable from $Y'=\mu _{\sigma (i)}(Y)$. From equation (3.16) and since $T_{pp', \sigma }(f_{p, x})$ is a sum of two monomials, then $P$ must be an integer. Finally, since the coefficients of $T_{pp', \sigma }(f_{p, x})$ and $f_{p', y_{\sigma (i)}}$ are all ones, then $P=1$.\\ Hence

 \begin{equation}\label{}
   T_{pp', \sigma }(f_{p, x})= f_{p', y_{\sigma (i)}}.
\end{equation}
Therefore
\begin{equation}\label{}
    T_{pp', \sigma}(\mu _{i}(x_{i}))=\mu _{\sigma (i)}(y_{\sigma (i)}), \ \text{for all} \  i\in [1, n].
\end{equation}

$(2) \Rightarrow (3)$ from Theorem 3.14 and $(3)\Rightarrow (1)$ is immediate.

\end{proof}
\end{thm}

\begin{cor}

\begin{enumerate} If $p$ and $p'$ are two labeled seeds in  $\mathcal{A}(X, Q)$ and $\sigma\in \mathfrak{S}_{n}$. Then
              \item $p$ and $p'$ are $\sigma$-similar if and only if $T_{pp', \sigma}$ is a cluster automorphism of $\mathcal{A}$.
              \item If $\mathcal{A}$ is  positive cluster algebra. Then, the following are equivalent
\begin{enumerate}
  \item $T_{pp', \sigma}$ is a $\mathcal{X}$-preserver;
  \item $p$ and $p'$ are $\sigma $-similar;
  \item $T_{pp', \sigma}$ is a cluster automorphism.
\end{enumerate}
            \end{enumerate}
\end{cor}
\begin{proof} Special cases of Theorems 3.14 and 3.15 by taking $\mathcal{B}(Y, Q')=\mathcal{A}(X, Q)$
\end{proof}

\subsection{The exchange group and  the group of cluster automorphisms}

Definition 3.8 and 3.10 give rise to the following two subgroups of  $Aut_{K}(\mathcal{F})$.

\begin{defns}
\begin{enumerate}
\item [(a)]
The exchange group, denoted by EAut$\mathcal{A}(X, Q)$, is the subgroup of $Aut_{K}(\mathcal{F})$ generated by  $\{T_{pp', \sigma}|  p, p' \in \text{Mut}(X, Q), \ \ \sigma \in\mathfrak{S}_{n}\}$.

\item [(b)][2] The  cluster automorphisms group  Aut$\mathcal{A}(X, Q)$ is the subgroup of $Aut_{K}(\mathcal{F})$  that consists of all cluster automorphisms of $\mathcal{A}(X, Q)$.
\end{enumerate}
\end{defns}

\begin{cor}
\begin{enumerate}
\item
\emph{Aut}$\mathcal{A}(X, Q)=\{ T_{pp', \sigma}| p, p' \ \text{are} \ \sigma-\text{similar  in} \ \mathcal{A}(X, Q), \sigma \in \mathfrak{S}_{n} \}$.

\item
 \emph{Aut}$\mathcal{A}(X, Q)=\text{\emph{EAut}}\mathcal{A}(X, Q)$ if and only if the $\sigma$-similarity relation on  \emph{Mut}$(X, Q)$ has only one equivalence class.
 \end{enumerate}
\end{cor}
\begin{exams}\begin{enumerate}
\item If $(X, Q)$ is a seed of rank $2$ then Aut$\mathcal{A}(X, Q)$=EAut$\mathcal{A}(X, Q)$.
\item Let
$Q=\xymatrix{
 \cdot_{1} \ar[d]_{(2,2)} \ar[d] \ar[r]^{(2,2)}&\cdot_{2} \\
  \cdot_{3}  \ar[ur]_{(2,2)}}$. Then  Aut$\mathcal{A}(X, Q)$=EAut$\mathcal{A}(X, Q)$.
\end{enumerate}
\end{exams}
In [2], the authors computed the cluster automorphism groups for cluster algebras of Dynkin and Euclidean types. Using Corollary 3.18 and part (1) of Example 3.19, we provide presentations for exchange groups and the cluster automorphisms groups of cluster algebras of  types $\mathbb{B}_{2}$ and $\mathbb{G}_{2}$.

\begin{exam}
\begin{enumerate}
\item
\textbf{The  group of cluster automorphisms} Aut$(X, \mathbb{B}_{2})$.

  \begin{equation}\label{}
   \nonumber \text{Aut}\mathcal{A}(X, \mathbb{B}_{2})=\{ T _{1}, T _{2} |\ \ T^{2} _{1}= T^{2} _{2}=1,  (T _{1} T
    _{2})^{3}=1\}.
    \end{equation}

\item
\textbf{The group of cluster automorphisms} $Aut(X, \mathbb{G}_{2})$.

\begin{equation}\label{}
    \nonumber \text{Aut}\mathcal{A}(X, \mathbb{G}_{2})=\{ T _{1}, T _{2} |\ \ T^{2} _{1}= T^{2} _{2}=1,  (T _{1} T
    _{2})^{4}=1\}.
\end{equation}
\end{enumerate}
\end{exam}

\end{document}